\documentclass[12pt, reqno]{amsart}
\usepackage{amsfonts,amssymb,graphicx,enumerate,ytableau,fourier}

\numberwithin{equation}{section}

\theoremstyle{definition}

\theoremstyle{plain}
\newtheorem{theorem}{Theorem}[section]

\newcommand{\CC}{\mathbf{C}}

\newcommand{\RR}{\mathbf{R}}

\newcommand{\ZZ}{\mathbf{Z}}

\newcommand{\OO}{\mathcal{O}}

\newcommand{\ttt}{\mathfrak{t}}

\begin{document}

\title[The diagonal coinvariant ring for type B]{Representation theory and the diagonal coinvariant ring of the type B Weyl group}

\author{Carlos Ajila}

\author{Stephen Griffeth}

\address{Instituto de Matem\'aticas \\
Universidad de Talca }

\begin{abstract}
We explain how to use representation theory to give a lower bound on the dimension of the quotient ring by type $B_n$ diagonal invariants that improves upon the current known lower bound $(2n+1)^n$ by a quadratic polynomial in $n$.
\end{abstract}

\thanks{We thank Iain Gordon for comments on a preliminary version of this paper. The second author acknowledges the financial support of Fondecyt Proyecto Regular 1190597.}

\maketitle

\section{Introduction}

\subsection{Diagonal coinvariant rings and the work of Haiman and Gordon} In \cite{Hai}, Mark Haiman stated a number of influential conjectures and observations having to do with the quotient $R_W$ of a polynomial ring $R=\CC[x_1,x_2,\dots,x_n,y_1,y_2,\dots,y_n]$ in two sets of $n$ variables each by the ideal generated by positive degree $W$-invariant polynomials, where $W=S_n$ is the symmetric group acting by simultaneous permutations of the variables. In particular, he conjectured (and later proved, in \cite{Hai2} and \cite{Hai3}) that $\mathrm{dim}(R_W)=(n+1)^{n-1}$. Since then, many mathematicians have contributed to the discovery and proof of increasingly refined theorems on the structure of $R_W$ as a bigraded $W$-module; for recent progress and references see \cite{vWi}.

As Haiman observed in \cite{Hai}, the analogous quotient ring $R_W$ for a general real reflection group $W$ is more mysterious: for instance, for $W=W(B_4)$ the dimension of $R_W$ is $9^4+1$, exactly one more than the naive expectation $9^4$. In general,  we do not have even a conjectural formula for its dimension. However, Haiman conjectured that the naive expectation gives a correct lower bound $\mathrm{dim}(R_W) \geq (h+1)^n$, where $h$ is the Coxeter number of $W$ and $n$ is its rank, and this conjecture was confirmed by Iain Gordon \cite{Gor} using representation theory (we note that \cite{GoGr} later proved the analogous theorem in the case of a complex reflection group). As far as we know, there have been no further improvements to this bound in the years since Gordon's work, and it seemed likely the extra dimensions appearing were not explained by the representation theory. Our object here is to report that in fact representation theory \emph{does} explain why the dimensions exceed our naive expectation, by using it to give an explicit lower bound for the error 
$$\epsilon(W)=\mathrm{dim}(R_W)-(h+1)^n$$ in case $W$ is a type B Weyl group, in particular proving that it is positive for all $n\geq4$. The Coxeter number for $W(B_n)$ is $2n$, and a consequence of our main result is that $$\epsilon(W(B_n))=\mathrm{dim}(R_W)-(2n+1)^n$$ grows at least like $n^2/4$. We have chosen the notation to indicate that the term $(h+1)^n$ should be regarded as the principal term in an approximation to the dimension of $R_W$; for a real reflection group $W$ we expect that $\epsilon(W)$ is very small compared to $(h+1)^n$. 

While the same ideas apply to all (complex) reflection groups $W$, we focus here on obtaining explicit lower bounds for $\epsilon(W)$ for $W=W(B_n)=G(2,1,n)$ the type B Weyl group, leaving more subtle calculations for future work. With this in mind, our main theorem is:

\begin{theorem} \label{main}
Let $W=W(B_n)$ be the group of signed permutations, acting on $$R=\CC[x_1,x_2,\dots,x_n,y_1,y_2,\dots,y_n]$$ by (simultaneous) sign-changes and permutations of the two sets of variables, let $I_W$ be the ideal of $R$ generated by the positive-degree $W$-invariant polynomials, and let $R_W=R/I_W$ be the quotient ring. Then
$$\epsilon(W(B_n)) \geq  \begin{cases} n(n-4)/4 \quad \hbox{if $n=0$ mod $4$,} \\ n(n-6)/4 \quad \hbox{if $n=2$ mod $4$, and} \\  (n-1)(n-3)/4 \quad \hbox{if $n$ is odd.} \end{cases}$$ Moreover $\epsilon(W(B_4)) \geq 1$ and $\epsilon(W(B_6)) \geq 3$.
\end{theorem}

\subsection{Our strategy} We now describe our strategy for a general complex reflection group $W$. There should be no surprises here for experts, up until we describe how to implement it and prove that it works. The strategy is this: each irreducible representation $L$ of the rational Cherednik algebra $H_c(W)$ with the property that the determinant appears exactly once in $L$ (we call such representations \emph{coinvariant type}) produces a lower bound for the graded dimension of each isotypic component of $R_W$ . Therefore the supremum of the graded characters over the (not linearly ordered!) set of coinvariant type representations of $H_c(W)$ is a lower bound for the graded character of $R_W$. We note here that it is not just \emph{one} coinvariant type representation giving complete information---rather, we should choose various $L$'s with different lowest weights and deformation parameters giving complementary information about different isotypic components of $R_W$. Gordon's choice of $L$ gives the \emph{principal term}, by far the largest contribution to the graded character of $R_W$. In this paper, the new information beyond Gordon's choice comes from the isotypic component for a certain linear character in a class of representations $L$ with lowest weights indexed by hooks for $n=5$ and $n>6$ (though other choices giving are possible in general and will  give additional information; this happens already for the cases $n=4$ and $n=6$, where the hooks do not provide information beyond Gordon's). The surprises (to us, at least) are: firstly, that there exist a great many coinvariant type representations giving new information, and secondly, that we now have enough control over the characters of coinvariant type representations to be able to say something new about $R_W$ using $L$'s that produce lower order terms in this way. Once one knows the area of search and has tools delicate enough to compute characters 
explicitly, the results are not hard to prove. This suggests the following program: classify the coinvariant-type representations (equivalently, the one-dimensional modules for the spherical rational Cherednik algebra), and for each, compute its character as explicitly as possible, thereby exhausting the information provided by our strategy. 

\section{Proof of Theorem \ref{main}}

\subsection{Outline} In this section we introduce notation and give the proof of Theorem \ref{main}, which is a corollary to the more detailed statement in Theorem \ref{mainprecise} below. We rely on \cite{Gri} and \cite{FGM} for technical tools, in particular identifying sufficiently many coinvariant type representations amongst the $\ttt$-diagonalizable representations classified in \cite{Gri}, and computing characters using the tools from \cite{FGM}.

\subsection{The rational Cherednik algebra} For a positive integer $n$ and deformation parameters $c,d \in \RR$, we will write $H_{c,d}$ for the rational Cherednik algebra of type $W(B_n)$, referring to \cite{Gri} (especially subsection 8.2) for the precise definition and conventions (which we believe to be totally standard). We note that $H_{c,d}$ is generated by two polynomial rings $\CC[x_1,\dots,x_n]$ and $\CC[y_1,\dots,y_n]$ and the group $W(B_n)$. A \emph{bi-partition} of $n$ is an ordered pair $\lambda=(\lambda^0,\lambda^1)$ of two partitions with $n$ total boxes; these index the isotypes of the irreducible complex representations $S^\lambda$ of $W(B_n)$ and of the irreducible objects $L_{c,d}(\lambda)$ of the category $\OO_{c,d}$ of $H_{c,d}$, among which occur all the representations of coinvariant type (as $c$ and $d$ vary). Given a representation $L=L_{c,d}(\lambda)$ of coinvariant type, we may filter $L$ by putting the determinant $\delta$ in degree $0$, and setting
$$L^{\leq m}=H_{c,d}^{\leq m} \cdot \delta,$$ where the filtration on $H_{c,d}$ is defined by demanding that the operators $x_i$ and $y_i$ have degree $1$ and each $w \in W$ has degree $0$. Then there is a surjection of graded $W$-modules (now with the Euler grading in which $x_i$ has degree $1$, $y_i$ has degree $-1$, and $W$ has degree $0$)
$$R_W \to \mathrm{gr}(L) \otimes \mathrm{det}^{-1}$$ of $R_W$ onto $\mathrm{gr}(L) \otimes \mathrm{det}^{-1}$ inducing an inequality (to be interpreted coefficient-wise) of graded $W$-characters
\begin{equation}
\mathrm{ch}(R_W) \geq \mathrm{ch}(L \otimes \mathrm{det}^{-1}),
\end{equation}
 where we have shifted the internal (Euler) grading on $L$ so that the determinant occurs in degree $0$ and as mentioned above we take the grading on $R_W$ in which $x_i$ has degree $1$ and $y_i$ has degree $-1$ for all $1 \leq i \leq n$.

\subsection{The detailed version of the main theorem}  For integers $1 \leq i \neq j \leq n$ we will write $(ij) \in W(B_n)$ for the transposition matrix interchanging coordinates $i$ and $j$ and $\zeta_i$ for the matrix changing the sign of the $i$th basis vector and leaving all other basis vectors fixed. Let $\chi$ be the linear character of $W(B_n)$ determined by $\chi((12))=-1$ and $\chi(\zeta_1)=1$ and let $\chi'$ be the product of $\chi$ and the determinant character $\mathrm{det}$ (so $\chi'((12))=1$ and $\chi'(\zeta_1)=-1$). We put
$$\epsilon_\chi(W(B_n))=\mathrm{dim}(R_W^{\chi'}))-\mathrm{dim}(L_{(2n+1)/2n,(2n+1)/2n}(\mathrm{triv})^\chi),$$ where for a character of $W$ we use superscripts to indicate istotypic components. The object of the remainder of the paper is to prove the next theorem, which by Gordon's result (c.f. \ref{Gordon} below) implies Theorem \ref{main}. In other words, the only gain we are exploiting so far comes from $\chi$; other representations may in general contribute to $\epsilon$ but the combinatorics of these is complicated enough to be left to future work.
\begin{theorem} \label{mainprecise}
For all $n \in \ZZ_{>0}$, 
$$\epsilon_\chi(W(B_n)) \geq  \begin{cases} n(n-4)/4 \quad \hbox{if $n=0$ mod $4$,} \\ n(n-6)/4 \quad \hbox{if $n=2$ mod $4$, and} \\  (n-1)(n-3)/4 \quad \hbox{if $n$ is odd.} \end{cases}$$ Moreover $\epsilon_\chi(W(B_4)) \geq 1$ and $\epsilon_\chi(W(B_6)) \geq 3$.
\end{theorem} 

\subsection{Tools} To prove Theorem \ref{mainprecise} we need to compute the dimension of the space of semi-invariants in specific examples of the modules $L_{c,d}(\lambda)$, for which purpose we rely on the classification theorem from \cite{Gri} and character formula from \cite{FGM}. We recall that $H_{c,d}$ contains a commutative subalgebra $\ttt$, the \emph{Dunkl-Opdam subalgebra}, generated by elements $$z_1,\dots,z_n,\zeta_1,\dots,\zeta_n,$$ and that \cite{Gri} Theorem 1.1 (a crucial step in the classification of the unitary representations) classifies those $L_c(\lambda)$'s on which $\ttt$ acts semisimply. Moreover, letting $H_{2,n}^{\mathrm{aff}}$ denote the subalgebra of $H_{c,d}$ generated by $\ttt$ and $W=G(2,1,n)$ (a certain generalization of the degenerate affine Hecke algebra), section 4 of \cite{FGM} describes the $H_{2,n}^{\mathrm{aff}}$ and $\CC W$-module structure of each $\ttt$-diagonalizable module explicitly. We will use only finite-dimensional $\ttt$-diagonalizable representations in what follows (fortunately, there are lots of these, and in fact we could have relied on the classification of the $\ttt$-diagonalizable standard modules in \cite{Gri2} for the exampes we need). We next describe the controlling combinatorics.

\subsection{Charged contents} Given a box $b$ of a partition, we define its \emph{content} to be 
$$\mathrm{ct}(b)=j-i$$ if it occurs in the $i$th row and $j$th column. We will define $d_0=d$ and $d_1=-d$, and in general $d_i=d_j$ for $i=j$ mod $2$ defines $d_i$ for all $i \in \ZZ$. For a box $b$ of a bi-partition $\lambda=(\lambda^0,\lambda^1)$ we write $\beta(b)=0$ if $b \in \lambda^0$ and $\beta(b)=1$ if $b \in \lambda^1$ (and more generally, define $\lambda^i$ for $i \in \ZZ$ by taking the residue of $i$ mod $2$). The \emph{charged content} of a box $b$ of a bi-partition $\lambda=(\lambda^0,\lambda^1)$ is then 
$$\mathrm{ct}_c(b)=d_{\beta(b)}+2\mathrm{ct}(b) c.$$ Given boxes $b,b' \in \lambda$, we write $b \leq b'$ if $\beta(b)=\beta(b')$ and $b$ is (weakly) above and to the left of $b'$.

\subsection{Tableaux} \label{tableaux} Given a bi-partition $\lambda$, we let $\mathrm{Tab}_{c,d}(\lambda)$ be the set of fillings $Q:\lambda \to \ZZ_{\geq 0}$ of the boxes of $\lambda$ by non-negative integers with the properties:
\begin{itemize}
\item[(a)] $Q$ is weakly increasing across rows and down columns in each component of $\lambda$,
\item[(b)] $Q(b) < k $ if $k$ is a positive integer with $$k=\mathrm{ct}_c(b)-d_{\beta(b)-k}, \quad \text{and}$$ 
\item[(c)] $Q(b) \leq Q(b')+k$ if $k$ is a positive integer with $k=\mathrm{ct}_c(b)-\mathrm{ct}_c(b') \pm 2 c$ for boxes $b, b'$ with $b' \in \lambda^{\beta(b)-k}$.
\end{itemize}  We refer to $Q$ as \emph{generic} if the inequalities which occur in (c) are strict. Given $Q \in \mathrm{Tab}_{c,d}(\lambda)$ we define another set $Q_{c,d}$ of tableaux to be the set of fillings $P$ of $\lambda$ with the properties
\begin{itemize}
\item[(a)] $P$ is a bijection from the boxes of $\lambda$ to $\{1,2,\dots,n\}$,
\item[(b)] if $b \leq b'$ and $Q(b)=Q(b')$ then $P(b)>P(b')$, and
\item[(c)] if $k$ is a positive integer with $k=\mathrm{ct}_c(b)-\mathrm{ct}_c(b') \pm 2 c$ for boxes $b, b'$ with $b' \in \lambda^{\beta(b)-k}$ and $Q(b)=Q(b')+k$ then $P(b)>P(b')$.
\end{itemize} Note that condition (c) is relevant only for non-generic $Q$. Theorem 1.1 of \cite{Gri} states that when $L_{c,d}(\lambda)$ is $\ttt$-diagonalizable, it has a basis $f_{P,Q}$ (these are the representation-valued Jack polynomials from \cite{Gri2}) indexed by pairs consisting of $Q \in \mathrm{Tab}_{c,d}(\lambda)$ and $P \in Q_{c,d}$, and in the proof of Theorem 4.1 of \cite{FGM} it is shown that in this case and assuming $c \neq 0$, the restriction of $L_{c,d}(\lambda)$ to the algebra $H_{2,n}^{\mathrm{aff}}$ is the direct sum of simple $H_{2,n}^{\mathrm{aff}}$-modules $L_Q$ for $Q \in \mathrm{Tab}_{c,d}(\lambda)$, with
$$L_Q=\CC \{f_{P,Q} \ | \ P \in Q_{c,d} \}.$$

\subsection{Simple $H_{2,n}^{\mathrm{aff}}$-module via skew diagrams} A \emph{diagram} is a finite subset $D \subseteq \CC^2$. We call $D$ \emph{skew} if for all $(a,b) \in D$ and all non-negative integers $s,t \in \ZZ_{\geq 0}$ such that $(a+s,b+t) \in D$, we have $(a+i,b+j) \in D$ for all $0 \leq i \leq s$ and $0 \leq j \leq t$. The \emph{connected components} of a diagram $D$ are the equivalence classes for the equivalence relation generated by $(a,b) \cong (a,b+1) \cong (a+1,b+1)$. As usual, we will visualize each connected component of a skew diagram as a collection of boxes, obtained (up to a translation) as the difference of two Young diagrams of (not unique) partitions. As for partitions, we define the \emph{content} of $(a,b) \in D$ by $\mathrm{ct}(a,b)=a-b$.

In section 3 of \cite{FGM}, it is explained that (for $c \neq 0$) the 
$\ttt$-diagonalizable simple $H_{2,n}^{\mathrm{aff}}$-modules are indexed by pairs of skew diagrams $D^0, D^1$ with $n$ total elements. We will abuse notation by referring to such a pair simply as a skew diagram when confusion will not result. A \emph{standard Young tableau} on a skew diagram $D=(D^0,D^1)$ of size $n$ is a bijection from the boxes of $D$ to the integers from $1$ to $n$ which is increasing across rows and down columns within each $D^i$. 

Fixing a simple $\ttt$-diagonalizable $H_{2,n}^{\mathrm{aff}}$-module $M$ and a non-zero vector $m\in M$ with 
$$z_i \cdot m=a_i m \quad \text{and} \quad \zeta_i \cdot m=(-1)^{b_i} m \quad \hbox{for some $a_i \in \CC$ and $b_i \in \{0,1 \}$,}$$ by subsection 4.5 and Theorem 3.2 of \cite{FGM}, there is a unique (up to diagonal slides of its connected components) skew diagram $D=(D^0,D^1)$ and standard Young tableau $T$ of shape $D$ such that $$\frac{a_i}{2c}=\mathrm{ct}(T^{-1}(i)) \quad \text{and} \quad b_i=\beta(T^{-1}(i)) \quad \hbox{for all $1 \leq i \leq n$.}$$ Here as above, $\beta(b)=i$ if $b \in D^i$. Moreover, $D$ depends only on the isotype of $M$, independent of the choice of $m$. We refer to the numbers $a_i$ and $b_i$ (with $b_i$ equal to $0$ or $1$) as the $i$th component of the \emph{$\ttt$-weight} of $m$.

\subsection{The structure of $L_Q$ as a $\CC W$-module} \label{LQ} Here we indicate (for $c \neq 0$, which will hold in all the examples we need) how one may use the preceding to compute the $\CC W$-module structure of $L_Q$ for each $Q \in \mathrm{Tab}_{c,d}(\lambda)$ (though in practice the combinatorics can be quite involved). For the proof of Theorem \ref{mainprecise} we will need only to know about occurrences of the linear characters $\mathrm{det}$ and $\chi$ in $L_Q$, which permits considerable simplification. 

Thus fix $Q \in \mathrm{Tab}_{c,d}(\lambda)$ and $P \in Q_{c,d}(\lambda)$. For $1 \leq i \leq n$, the $i$th component of the $\ttt$-weight of $f_{P,Q}$ is $(a_i,b_i)$ with
\begin{equation} \label{weight sequence} a_i=Q(P^{-1}(i))+1-(d_{\beta(P^{-1}(i))}-d_{\beta(P^{-1}(i))-Q(P^{-1}(i))-1})-2 \mathrm{ct}(P^{-1}(i)) c\end{equation} and $b_i=\beta(P^{-1}(i))-Q(P^{-1}(i))$. We write $D_{Q,c,d}=(D_{Q,c,d}^0,D_{Q,c,d}^1)$ for the skew diagram obtained from this as in the previous subsection (it is independent of the choice of $P \in Q_{c,d}$). By 3.10 of \cite{FGM} the occurrences of the $W$-module indexed by a bipartition $(\mu^0,\mu^1)$ in $L_Q$ are indexed by pairs $T^1,T^2$ of Littlewood-Richardson tableaux on $D^0_{Q,c,d},D^1_{Q,c,d}$ with weights $\mu^0$ and $\mu^1$. 

\subsection{The special case of linear characters} \label{linear case} Combining the Littlewood-Richardson rule with the preceding shows that $\mathrm{det}$ (resp., $\chi$) occurs in $L_Q$ precisely if $Q(b)=\beta(b)+1$ (resp., $Q(b)=\beta(b)$) for all $b \in \lambda$ and moreover no two boxes of $D_{Q,c,d}^1$ (resp., $D_{Q,c,d}^0$) are in the same row (that is, the diagram is a \emph{vertical strip}), in which case it appears with multiplicity one. This can happen only if $Q$ is row-strict (strictly increasing across rows); if moreover $Q$ is generic then row-strictness is sufficient (this follows e.g. from Theorem 2.2 of \cite{DuGr}). We will use only these consequences of \ref{LQ} from now on.

\subsection{First example: the Gordon module} \label{Gordon} To illustrate the previous constructions and for later use, we now take $\lambda=((n),\emptyset)$ to be the bi-partition indexing the trivial representation, and choose parameters $c$ and $d$ generic subject to $$2n+1=2(d+(n-1)c)$$ (Gordon's choice $c=d=(2n+1)/2n$ is generic for our purposes, as follows from \cite{Gri} Theorem 1.1). With such a choice of parameters $L_{c,d}(\mathrm{triv})$ is $\ttt$-diagonalizable, and letting $b$ be the unique removable box of $\lambda$, we have
$$\mathrm{Tab}_{c,d}(\lambda)=\left\{Q \ \text{weakly increasing} \ | \ Q(b) \leq 2n \right\}.$$ All such $Q$'s are generic; by \ref{linear case} the determinant occurs precisely for $Q$ with strictly increasing odd entries $1,3,\dots,2n-1$, so that $L_{c,d}(\mathrm{triv})$ is of coinvariant type, and $\chi$ occurs exactly for $Q$ with strictly increasing even entries chosen from $0,2,4,\dots,2n-2,2n$. There are $n+1$ such $Q$ (corresponding to a single $\mathfrak{sl}_2$-string). 
 
\subsection{The case $\lambda=((2,2),\emptyset))$.} Here we explain that for $n=4$ the extra dimension in $R_W$ comes from $L_{c,5/2}(\lambda)$ for $c$ generic. The equation $5=2(d+\mathrm{ct}(b) c)$ holds, where $b$ is the unique removable box, and imposes $Q(b) \leq 4$ for $Q \in \mathrm{Tab}_{c,d}(\lambda)$ (all are generic); the representation is of coinvariant type since the unique occurrence of the determinant comes from  
$$Q= \begin{ytableau}
1 & 3 \\ 1 & 3 
\end{ytableau}$$
and there are six occurrences of $\chi$, corresponding to the following fillings $Q$:
$$\begin{ytableau}
0 & 2 \\ 0 & 2 
\end{ytableau} \quad \begin{ytableau}
0 & 2 \\ 0 & 4 
\end{ytableau} \quad \begin{ytableau}
0 & 4 \\ 0 & 4 
\end{ytableau} \quad \begin{ytableau}
0 & 2 \\ 2 & 4 
\end{ytableau} \quad 
\begin{ytableau}
0 & 4 \\ 2 & 4 
\end{ytableau} \quad
\begin{ytableau}
2 & 4 \\ 2 & 4 
\end{ytableau}$$ When combined with \ref{Gordon} this shows $\epsilon_\chi(W(B_4)) \geq 1$. A similar analysis with $\lambda=((3,3),\emptyset)$ and generic parameters subject to $7=2(d+c)$ shows  $\epsilon_\chi(W(B_6)) \geq 3$. 

\subsection{Hook lowest weights} Next we consider the hook bipartition $\lambda=((k,1^m),\emptyset)$ for positive integers $k$ and $m$ with $k+m=n$. For the first time we must use uniquely determined parameters (the modules giving us the information we need live only at isolated points in the parameter space, unlike in the previous cases). Now $\lambda$ has two removable boxes; let $b$ be its removable box of larger content and let $b'$ be its removable box of smaller content. We choose the deformation parameters $c$ and $d$ so that conditions (b) and (c) from the definition of $\mathrm{Tab}_{c,d}(\lambda)$ in \ref{tableaux} are quite restrictive:
$$3=2d+2 \mathrm{ct}(b') c \quad \text{and} \quad 2 n c=2k.$$ Then a weakly increasing filling $Q$ of $\lambda$ belongs to $\mathrm{Tab}_{c,d}(\lambda)$ if and only if $Q(b') \leq 2$ and $Q(b) \leq Q(b')+2k$, and $Q$ is generic if $Q(b) \leq Q(b')+2k-1$.

In order that $c$ have denominator precisely $n$ (which implies, by Theorem 1.1 of \cite{Gri}, that $L_{c,d}(\lambda)$ is $\ttt$-diagonalizable), we must choose $k$ prime to $n$. Moreover, in order to maximize the quantity $2+km$ which appears below, we'd like $k$ to be roughly $n/2$. Here is one way to achieve both aims:
$$k=\begin{cases}  n/2+1 \quad \hbox{if $n=0$ mod $4$,} \\ n/2+2 \quad \hbox{if $n=2$ mod $4$, and} \\ (n+1)/2 \quad \hbox{if $n$ is odd.}
\end{cases}$$ We now establish that with these choices $L_{c,d}(\lambda)$ is of coinvariant type.

There is a unique generic row-strict $Q$ with odd entries: it has $1$'s in the column and strictly increasing odd entries $1,3,5,\dots,2k-1$ across the row. For instance for $k=3$ and $m=2$ it is
$$\begin{ytableau}
1 & 3 & 5 \\ 1 \\ 1
\end{ytableau}.$$ This produces a copy of the determinant. But one observes that the other row-strict $Q$'s with odd entries  in $\mathrm{Tab}_{c,d}(\lambda)$ have $Q(b)=2k+1=Q(b')+2k$, and  using \eqref{weight sequence} (it is convenient to choose $P \in Q_{c,d}$ with $P(b)=n$, $P(b')=1$) together with \ref{linear case} implies that the diagram $D_{Q,c,d}^1$ has two boxes in the same row and therefore does not produce a copy of the determinant. Hence $L_{c,d}(\lambda)$ is of coinvariant type. 

To give a lower bound for the number of occurrences of $\chi$ we observe that row-strict tableaux $Q$ of shape $(k,1^m)$ with even entries taken from $\{0,2,\dots,2k\}$ satisfying $Q(b') \leq 2$ and $Q(b) \leq Q(b')+2k-2$ are all generic, so produce copies of $\chi$. One checks that there are $2+km$ of these. But as we remarked in \ref{Gordon}, in Gordon's representation $L_{(2n+1)/2n,(2n+1)/2n}(\mathrm{triv})$, the isotype $\chi$ appears only $n+1$ times. The difference $2+km-(n+1)$ is therefore a lower bound for the error term $\epsilon_\chi(W(B_n))$, which proves Theorem \ref{mainprecise} and consequently Theorem \ref{main}. We remark that using only rectangular partitions as in the cases $n=4$ and $6$ would have given (for even $n$) a bound that grows like $n^2/8$. 

\def\cprime{$'$} \def\cprime{$'$}

\end{document}